\documentclass[12pt]{article}

\usepackage[utf8]{inputenc}
\usepackage[T1]{fontenc}
\usepackage{amsmath,amssymb}
\usepackage{url}
\usepackage{xcolor}

\title{Hostages}
\author{Ernesto Mordecki\footnote{Centro de Matemática, Facultad de Ciencias, Uruguay}}
\begin{document}
\maketitle
\begin{abstract}
Our house has been taken over\footnote{An allusion to Julio Cortázar's short
	story ``Casa tomada''~\cite{Cortazar1951}.}.
It is not completely destroyed, but many of its essential functions and
facilities no longer work as they used to. Essentially, they are out of our
control. Artificial intelligence, through its large language models (LLMs), has
taken it.
And we, the inhabitants of the house,
without really noticing what was happening, and surprised by the effectiveness of LLMs,
have handed over our most sacred belongings: the capacity to think
about a problem, the joy of having an idea that solves it, the
beautiful effort needed to consolidate an idea in our common language.
And there is no turning back.

We are mathematicians, and our whole world is in check.
We wonder: is it checkmate? Can we take back control of our house, which will
undoubtedly not work in the same way?

Paradoxically, the weapons that attacked us are built on some of our own ideas.
Of course, it is not entirely our work.
The computer science knowledge needed to develop these large models evolved
in both its theoretical and its material aspects.
As technology evolved (infinitely small processors, infinitely large computing centers),
turning quantity into quality,
computer scientists encoded into software the ideas needed to produce these tools
(let us not call them weapons).

These notes look into the past and the future of the mathematical profession.
What will inevitably happen?
What will remain and what will disappear?
And, in consequence,
what is our line of action?
To organize the ideas, after some introductory historical comments, we try to answer the questions:
(1) Is mathematics necessary for the development of our society?
(2) Will mathematicians be necessary in the near future?
(3) How the profession will change, through a list of \emph{Exercises}: open questions and possible actions.
\end{abstract}

\section{Background}


We begin with a brief discussion of some selected technological advances in human history,
ending with a comparison between the human brain and artificial language models.

\subsection{Writing}

In ancient times, oral tradition was the way to transmit knowledge, which was stored in human brains.
With the emergence of writing, knowledge was no longer stored only in the human
brain, but began to be stored on physical supports: handwritten books.
This was an important moment in the
development of humankind;
it is considered the passage from prehistory to history of a civilization.
This process was not a simple one, and, of course, there was resistance to it.
A classical example is the point of view of Socrates, one of the greatest philosophers of all time, who was against writing and defended oral tradition (see \cite{wiki1}).
Paradoxically again, the ideas of Socrates have been preserved thanks to the writings of Plato, his pupil.

\subsection{The printing press}

The reproduction of books had been the work of scribes since Antiquity and throughout the Middle Ages,
up to the 15th century.
Although today we may think of this labor as boring and repetitive,
scribes were very learned people for their time:
they were able to read and write, and they knew and understood the scarce books that reached them \cite{Vallejo2019}.
Gutenberg's invention of the printing press around 1455 made it possible to produce hundreds of identical copies
at a fraction of the cost. 
Within a few decades, the manual copying of books ceased to be the main means
of transmitting knowledge, and the scribes' trade began to disappear~\cite{Vallejo2019}.

\subsection{Chess machines}
Jumping to the last century, Garry Kasparov, then world chess champion, played a match against
Deep Blue\footnote{Deep Blue was a chess-playing computer built by IBM.}
in May 1997.
Deep Blue defeated Kasparov $3\tfrac{1}{2}$--$2\tfrac{1}{2}$ in a six-game
match,
the first time a computer beat a reigning world champion under tournament conditions.

The result was a sensation. Once a computer could beat the most
talented human, chess lost its aura as one of the most relevant
intellectual activities of humankind and became something different:
a simple game, an interesting game, a deep
game, but nothing more than a game.\footnote{
Similar examples can be found in checkers, with the program Chinook in the 1990s,
and in Go, with AlphaGo in 2015--2016.}

\subsection{The brain and neural networks}\label{sec:brain}

To assess whether machines can substitute for human thinking,
we make some crude comparisons.
The human brain contains about $8.6\times 10^{10}$
neurons~\cite{Azevedo2009} and on the order of $10^{14}$--$10^{15}$
synapses, and it operates on roughly $20$ watts.
The largest current artificial networks have on the order of $10^{11}$--$10^{12}$ parameters,
which are often compared to synapses.
They are trained on the order of
$10^{13}$ words or tokens, consuming megawatts over months. A human, by
contrast, encounters perhaps $10^{8}$--$10^{9}$ words by adulthood.
Although they are very different objects, when regarded as machines
that receive, store, process and output information, they show clear
similarities, at least in terms of their input--output behavior.
And, undoubtedly, artificial neural networks are expected to keep growing
in size and capacity, to limits that are difficult to assess today.
It is now hardly disputed that machines can pass the Turing test \cite{wiki2}, and we have to recognize that
machines show a certain ability to exhibit ``intelligent'' behavior, in the undefined sense considered by Turing.

\section{Is mathematics necessary?}

The main question we want to answer in these notes is whether mathematicians will be necessary in the near future
and, if the answer is positive in some sense, how mathematical activity will develop.
The first of these questions will be answered by our employers,
not by ourselves.
Nevertheless, to understand how our profession is changing, to try to figure out
what it will look like in the near future,
and especially to understand how to act now, it is necessary to examine it.

We therefore begin with a more elementary question, and an easier one to answer: is mathematics necessary (for our
society, for our lives)?
And, if the answer is positive,
can we, as a society,
preserve and enhance mathematics without mathematicians?

The presence of mathematics in society and everyday life is ubiquitous,
and it is difficult to imagine how our collective life
could continue as we understand it without mathematics.
Let us examine its different roles in society.

\subsection{Technology}

First, mathematics is a determining ingredient of \emph{technology and production}:
behind practically any object or service there is some hidden mathematical structure or achievement.
Cars, planes, ships, navigation instruments, clocks, phones, etc.
If not impossible, it seems difficult to find objects or services free of mathematics (not to mention computers, the internet, cell phones, and a long list of other achievements of the 20th century).

Artificial intelligence itself is perhaps one of the most spectacular achievements
with a mathematical basis of the 21st century.
Large language models (LLMs), the ones that answer our mathematical questions,
write and run our computer programs,
and prove our theorems,
have mathematical bricks as their central components, some of them relatively elementary:
neural networks are compositions
of linear and nonlinear multidimensional functions (as taught in advanced university calculus);
they are trained using the gradient descent algorithm, based on the chain rule (which we teach in the same advanced course);
they represent words as vectors and measure the distances between them (linear algebra);
and they compute the most probable word to follow a given one (probability and statistics).

One very relevant remark should be made, and it is not mathematical: the success of these LLMs
is also based on their large size, astonishingly large compared with what mathematicians and statisticians were used to
considering in their examples. This also contradicts some traditional modeling principles of statistics,
such as the principle that one should model with as few parameters as possible. This principle (implemented
by the Akaike information criterion, the Bayesian information criterion and subsequent improvements) is no longer applied.
It is also relevant to mention that many of the spectacular results of LLMs are the outcome of empirical work (for instance, the design of the network architecture), have no mathematical proof,
and their quality is, for now, difficult to assess or to compare objectively (for instance, in artificial image generation).

\subsection{Culture}

Second, mathematics is part of \emph{science and culture}.
Besides being part of human knowledge in itself, mathematics is an
indispensable tool in many related scientific areas; as Galileo Galilei poetically stated,
``[the universe] $\dots$ is written in the language of mathematics, and its characters are triangles, circles, and other geometric figures'' \cite{wiki3}.
It is also one of the oldest continuous intellectual traditions of humanity, comparable to
literature or music.


\subsection{Curriculum}

Third, mathematics is a \emph{curricular backbone} at all levels of our educational system.
Leaving aside primary and secondary school (somewhat beyond my reach),
mathematics is currently a central topic in the curricula of engineering,
physics, astronomy, chemistry, biology, economics, computer science and so on.
Almost every scientific and technological degree program
in the world begins with calculus and linear algebra, and later adds
differential equations, probability and statistics, at different levels of depth depending
on the prior preparation of the students and the expectations for their future development.

One question we have to face is the following: could a professional, with the
assistance of an AI system, perform their tasks, with enough knowledge of
their discipline to detect wrong or misplaced results, without
the need to learn how to compute areas and volumes,
the rules of differentiation, and so on?
It is then important to understand the process of learning mathematics in this new environment.

\subsection{Research}

Fourth, mathematics is a
\emph{research activity}.
As Thurston characterized it, ``mathematics has an essential recursive quality'' \cite{Thurston1994}.
We interpret this statement in the following way: to apprehend mathematics, to know and understand it deeply, to be able to improve and transmit human mathematical knowledge,
it is necessary to do research, to formulate new questions, and to try to answer them.
This is the activity that mathematicians carry out, individually and collectively,
and (as Thurston also illustrates) one of its multiple outcomes is the publication of a paper (see the interesting
metaphor with goals in soccer \cite{Thurston1994}).
In conclusion, the central role of mathematicians
is to know, preserve, reorganize, enlarge and transmit a body of knowledge across
generations. In this process, they organize seminars and conferences, publish new results in the form of articles,
synthesize new or old theories through the publication of books, and carry out other scientific and outreach activities.
Consequently, mathematical research is necessary in order to have people who understand mathematics,
who are able to transmit it to new generations,
and who also serve as a reference for applications and interactions with other disciplines and technologies.

\section{Will mathematicians be necessary in the near future?}

This second question is more difficult.
We explained above that, to transmit mathematical knowledge as we understand it today,
human mathematical research is necessary.
Also, to maintain the systems that
substitute part of our work (LLMs) and to make them more powerful, mathematicians also seem to be necessary.
Nowadays, the big AI companies have large research departments that include mathematicians.
But it is clear that this is not the traditional work of a mathematician at a university or a research institute.

At the moment, it is not clear whether LLMs can improve their
performance without human intervention.
But it may be the case that they can improve their performance with the
participation of technicians, engineers and computer scientists, and perhaps
that will be enough.

Regarding society in general,
in the old days (yesterday, last year, before the emergence of LLMs)\footnote{For me, the cutoff date was April 1st, 2026, when I subscribed to a Claude plan.}
mathematicians played a significant role in the transmission of mathematical knowledge to many future professionals,
such as engineers, physicists, economists, accountants, biologists, geologists, and many others, at the university level,
and in training new generations of mathematicians, to ensure the continuation of this work.

As discussed above,
it is important to discover and understand the new ways of mathematical thinking and learning,
in order to update the teaching proposals in strongly math-related degree programs.
A similar situation can be encountered in practically all professions (excluding mathematicians).

\section{Exercises}

It is, of course, premature to draw firm conclusions from the above comments.
But we are witnessing a fast-moving scenario,
without an opportunity to stop the process and think collectively
about the problems posed by AI to mathematics. The moment is now.
Trying to think about what our profession will look like in the near future,
we formulate an incomplete list of questions and statements to face it.

\begin{itemize}
\item\textbf{Quality over quantity.} Now that the production of articles and books has become extremely easy compared with yesterday (without AI), we should \textbf{definitely} return to assessing quality over quantity.
This applies both to the production of our mathematical results and to the selection of master's or PhD candidates,
or of candidates for new positions.
\item\textbf{Incorporating AI.} If we decide to incorporate AI into our work, we have to learn how to use it
and, more difficult, to discover new ways of working and learning mathematics.
At first glance, it seems that we have been relegated to reading and digesting what an LLM dictates about our problems.
This fact has two faces: on the one hand, it puts us in a secondary position; the ideas are no longer exclusively ours.
On the other, it gives us the possibility to solve traditional problems faster, and also to consider more difficult and central questions, out of our reach without AI.
Finally, it will result in an interaction that we should learn how to lead.

\item\textbf{Solving mathematical problems.} The traditional production of new results and solution of conjectures, i.e., the ``problem-solving'' activity of mathematicians, will change drastically. The production of new theories or the formulation of essentially new problems does not seem to be within reach of AI systems yet, but it is probably a matter of time.
\item\textbf{Publication of articles in refereed scientific journals.}
The traditional system, consisting of the publication of peer-reviewed articles to build a mathematical career, does not seem to be feasible in the AI era.\footnote{Currently, some journals have temporarily suspended the submission of papers ``in response to a substantial increase in manuscript submissions across the discipline, driven largely by the rapid development of artificial intelligence'' \cite{esaim}.}\textsuperscript{,}\footnote{It is not feasible to referee AI-written papers without AI, as some editorial boards recommend.}
How can we assess quality, acknowledge authorship, and preserve and disseminate mathematical knowledge in a new way?
How can we assess mathematical ability without the traditional ``list of publications'' in CVs?
\item\textbf{Preprints.} Preprints (such as those on arXiv or HAL) become more relevant.
As the production of results speeds up, and peer review (if possible at all) remains at the same pace,
the prepublication of results gains a new role. It is then important to make the results obtained publicly available
(of course, in careful human language) in repositories accessible to the community, and if possible in more than one place.
\item\textbf{Postgraduate studies.} Not only undergraduate studies will change course.
How should a PhD thesis in mathematics be supervised? What will a PhD dissertation and its defense look like?
\item\textbf{Mathematics of AI in the curriculum.}
As with any developing technology, AI itself generates a need for mathematical work and ideas.
In order to understand the general principles of the area, if not to participate in finding new results and producing new methods, we propose including the following contents in the curricula (for pure and applied mathematicians):
\begin{itemize}
\item Large language models and other mathematical methods of AI.
\item Automatic proof verification systems: the study and development of classical abstract areas of computer science and the formalization of mathematics.
\end{itemize}
\item\textbf{Hardware.} Installation, in mathematics departments, of autonomous LLMs that can be trained with local information, independent of commercial models.
\item\textbf{Research topics in the mathematics of AI.}
Part of a comprehensive response to the current changes is to involve the mathematical community in the production of AI methods, the proof (or disproof) of how they actually work, the development of new concepts in high-dimensional statistics to cope with these large models, and many other possible topics.
\end{itemize}

\end{document}